%% file: FCAA2013_Optimal_Random_Search.tex
 \documentclass[twoside,reqno,11pt]{fcaa}
 \input fcaa_style

 \usepackage{hyperref} 
 \usepackage{multirow}
 \usepackage{slashbox}
 \def\theequation{\arabic{section}.\arabic{equation}}

 \def\themycopyrightfootnote{\vspace*{3pt}
 \copyright \, Year\,  Diogenes  Co., Sofia
   \par  \noindent pp. xxx--xxx
   \hfill  \vspace*{-36pt}
  }

 \title[Optimal random search, fractional dynamics~$\cdots$]{Optimal random search, fractional dynamics and fractional calculus}
 \author[C. Zeng \& Y.Q. Chen]{Caibin Zeng$^{1}$, YangQuan Chen$^2$}

\begin{document}

 \vbox to 2.5cm { \vfill }


 \bigskip \medskip

 \begin{abstract}

What is the most efficient search strategy for the random located
target sites subject to the physical and biological constraints?
Previous results suggested the L\'evy flight is the best option to
characterize this optimal problem, however, which ignores the
understanding and learning abilities of the searcher agents. In the
paper we propose the Continuous Time Random Walk (CTRW) optimal
search framework and find the optimum for both of search length's
and waiting time's distributions. Based on fractional calculus
technique, we further derive its master equation to show the
mechanism of such complex fractional dynamics. Numerous simulations
are provided to illustrate the non-destructive and destructive
cases.
 \medskip


 \smallskip

 \end{abstract}

 \maketitle

 \vspace*{-16pt}


\section{Introduction}\label{Sec:1}

Over recent years the accumulating experimental evidences show that
the moving organisms are ubiquitous. For instance, the foraging
behavior of the wandering albatross (Diomedea exulans) on the ocean
surface was found to obey a power-law distribution
\cite{Viswanathan1996}; the foraging patterns of free-ranging spider
monkeys (Ateles geoffroyi) in the forests of the Yucatan Peninsula
was also found to be a power law tailed distribution of steps
consistent with a L\'{e}vy walk \cite{Ramos2004,Boyer2004}. More
experimental findings can be found in \cite[Part
II]{Viswanathan2011}. A number of foundational but important
questions arise naturally: How to model these organisms¡¯s movement
trajectories? What factors determine the shape and the statistical
properties of such trajectories? How to optimize the efficiency to
search of randomly located targets? These questions have been
studied in many field. For example, L\'{e}vy flight search was
claimed to be an optimal strategy in sparsely target site with an
inverse square power-law distribution of flight lengths
\cite{Viswanathan1999}. Then composite Brownian walk searches were
found to be more efficient than any L\'{e}vy flight when searching
is non-destructive and when the L\'{e}vy walks are not responsive to
conditions found in the search \cite{Benhamou2007}. In particular,
the movement patterns have scale-free and super-diffusive
characteristics. So the fractional Brownian motions and fractional
L\'{e}vy motions are possible to account for the movement patterns
\cite{Reynolds2009}. However, the above strategies ignore one
important factor, waiting time between the successive movement
steps, since the search agents need some time to understand the
visited target sites. The relation between waiting time and flight
length for efficient search was discussed \cite{Koyama2008}.

This fact inspires us to propose a potentially transformative
framework for optimal random search based on continuous time random
walk (CTRW). The CTRW strategy is composed of the flight lengths of
a movement step with a random direction, as well as the waiting time
elapsing between two successive movement steps, both of which are
independent random variables, identically, distributed according to
their probability densities. In addition, CTRW is also the
stochastic solution of non-integer order diffusion equation based on
fractional calculus \cite{Weiss1994}, which is a part of mathematics
dealing with derivatives and integration of arbitrary order
\cite{Podlubny1999}. Different from the analytical results of linear
integer-order differential equations, which are represented by the
combination of exponential functions, the analytical results of the
linear fractional order differential equations are represented by
the Mittag-Leffler function, which exhibits a power-law asymptotic
behavior \cite{Monje2010}. Therefore, fractional calculus is being
widely used to analyze the random signals with power-law
distributions or power-law decay of correlations \cite{Sheng2012}.
By choosing the flight lengths subject to heavy-tailed distribution
and finite characteristic waiting time, CTRW compasses L\'{e}vy
flight as a special case. Therefore, in this paper we propose the
CTRW optimal search framework, which may provide new insights into
the optimal random search in unpredictable environments.

The paper is organized of as follows: In Section \ref{sec:2}, we
review the L\'evy flight optimal random search strategy. Then we
formulate the CTRW optimal random search framework and find the
optimum for both of search length's and waiting time's distributions
in Section \ref{sec:3}. Based on fractional calculus, we derive the
corresponding master equation to CTRW search strategy in Section
\ref{sec:4}. Finally, we give some concluding remarks and close this
paper in Section \ref{sec:5}.

\section{L\'evy flight optimal random search}
\label{sec:2}

In this section, we review the basic idea of L\'{e}vy flight optimal
random search and reproduce the main results, which can help us
better understand the proposed CTRW strategy. In
\cite{Viswanathan1999,Viswanathan2000,Viswanathan2001}, the authors
assumed the search length distribution
\begin{equation}
    p(l_j)\sim l_j^{~-\mu},
\end{equation}
with $1<\mu<3$. Then they defined the search efficiency function
$\eta(\mu)$ to be the ratio of the number of target sites visited to
the total distance traversed by the forager as following
\begin{equation}
\eta(\mu)=\frac{1}{N\langle L\rangle},
\end{equation}
in which $N$ denotes the mean number of flights taken by a L\'{e}vy
forager in order to travel between two successive target sites and
$\langle L\rangle$ denotes the mean flight distance. For the case of
destructive foraging, they found that the mean efficiency $\mu$ has
no maximum, with lower values of $\mu$ leading to more efficient
foraging. For the case of nondestructive foraging, they found that
an optimal strategy for a forager is to choose $\mu^\ast=2$ when
$\lambda$ is large but not exactly known. The main results were
reproduced in the cases of destructive foraging and nondestructive
foraging, shown in Figure 2.1.

 \begin{center}
\includegraphics[height=2in,width=2.4in]{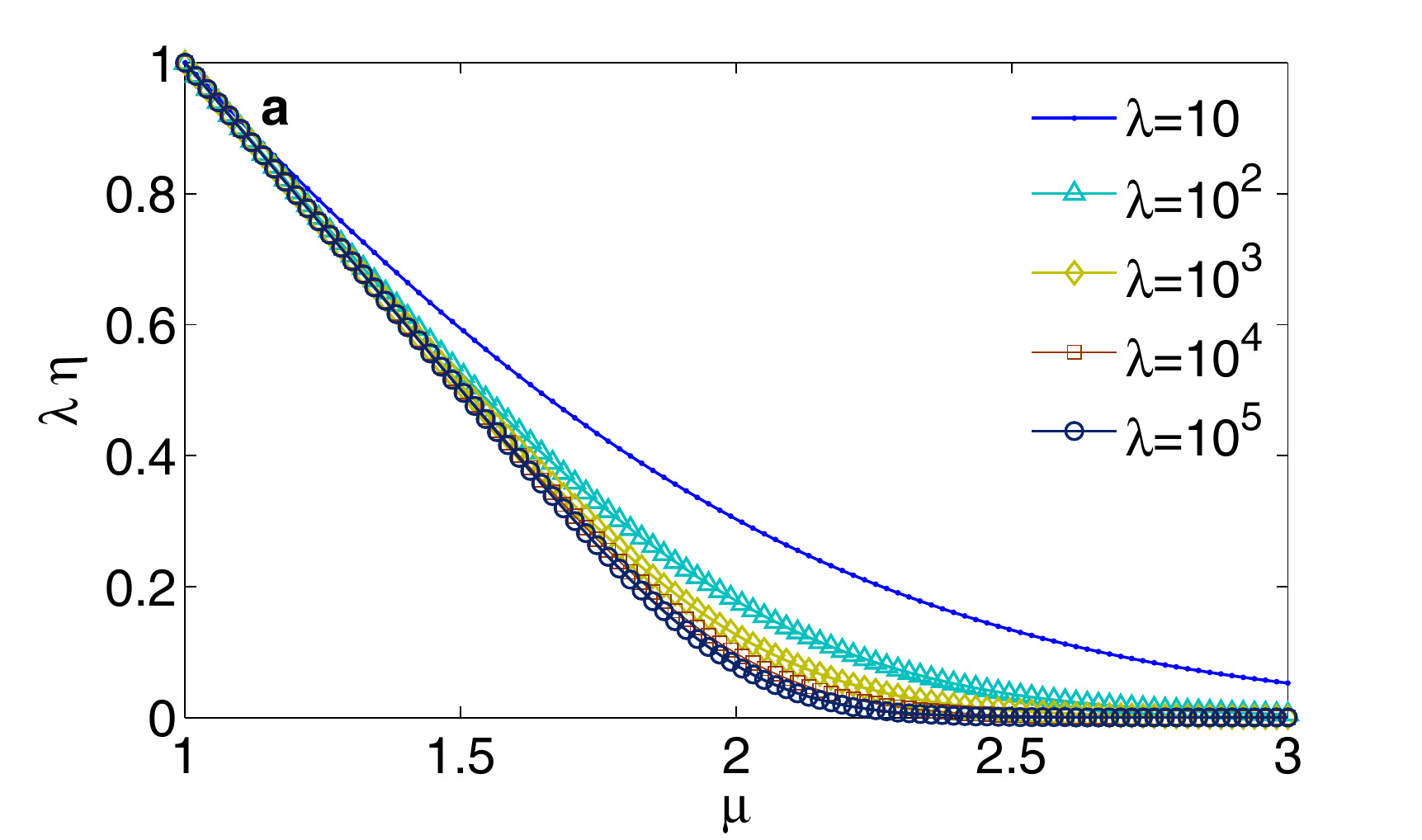}
\includegraphics[height=2in,width=2.4in]{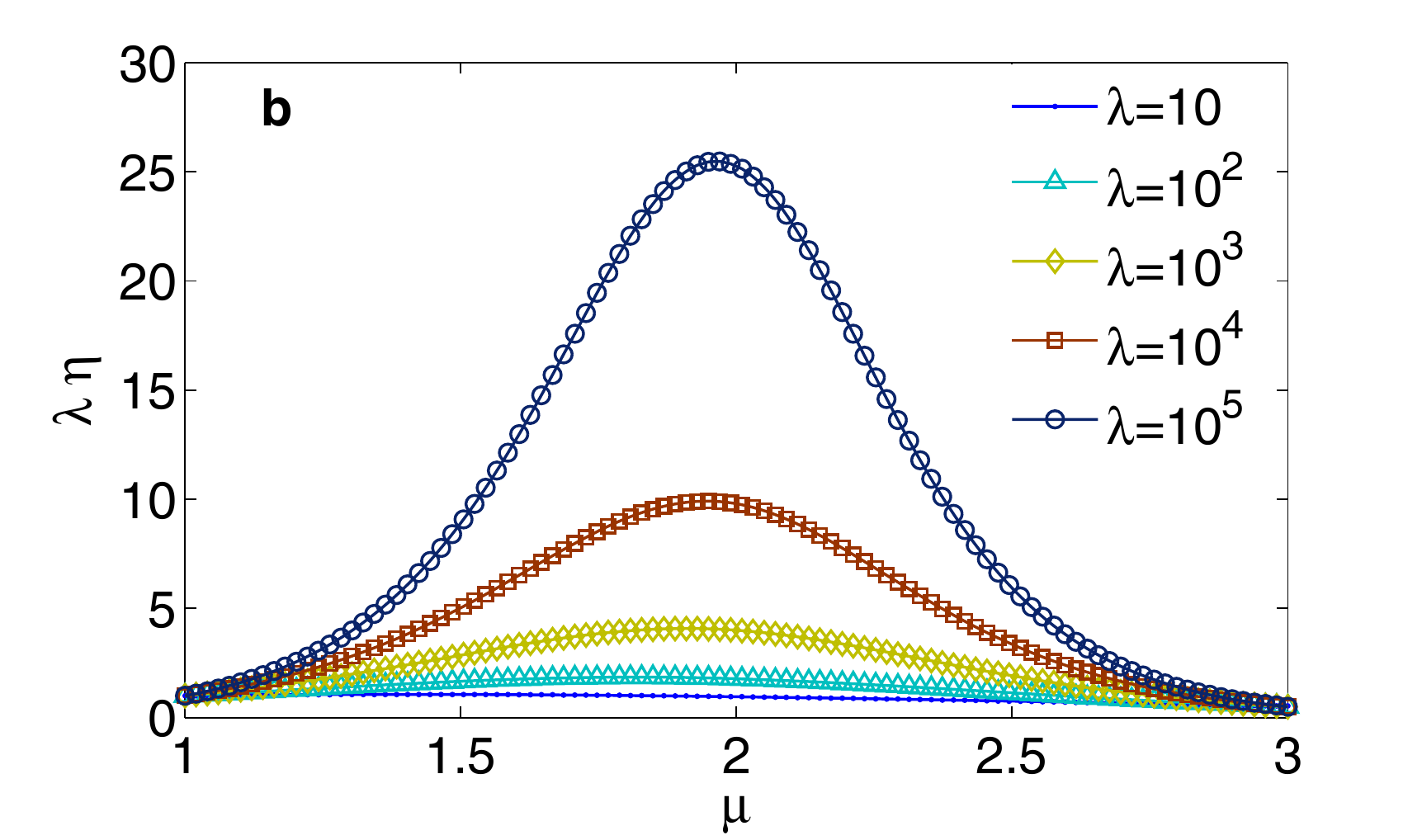}
 \bigskip

Fig. 2.1: The product of the search efficiency $\eta$ and the mean
free path $\lambda$ VS the parameter $\mu$ for different $\lambda$:
(a) the destructive case; (b) the nondestructive case.
 \end{center}

\section{CTRW optimal random search}\label{sec:3}

Now we state the main idea of CTRW search strategy by using the
probability distribution functions of search length and waiting
time. Specifically, the former is given by
\begin{equation}\label{eq2.1}
    w(l_j)\sim l_j^{~-(\alpha+1)},
\end{equation}
where $0<\alpha<2$, and $l_j$ is the the search length at the $j$th
step. The latter is characterized by
\begin{equation}\label{eq2.2}
    \psi(t_j)\sim t_j^{-(\beta+1)},
\end{equation}
where $0<\beta<1$, and $t_j$ is the waiting time length before
starting the $j$th step. Then the CTRW search strategy is described
by the following two simple rules: 1) If a target site is located
within a visible and finite distance $r_v$, then the search agent
moves on a straight line to it without learning; (2) if there is no
target site within a finite distance $r_v$, then the agent spends
some waiting time $t_j$, which is also can be characterized by a
power-law function (\ref{eq2.2}), to understand what is detected,
and chooses a random direction and a random distance $l_j$ following
another power-law distribution function (\ref{eq2.1}); otherwise, it
proceeds to the target as in first rule.

The schematic idea of CTRW search strategy is drawn in Figure 3.1.

 \begin{center}
\includegraphics[height=2in,width=2.4in]{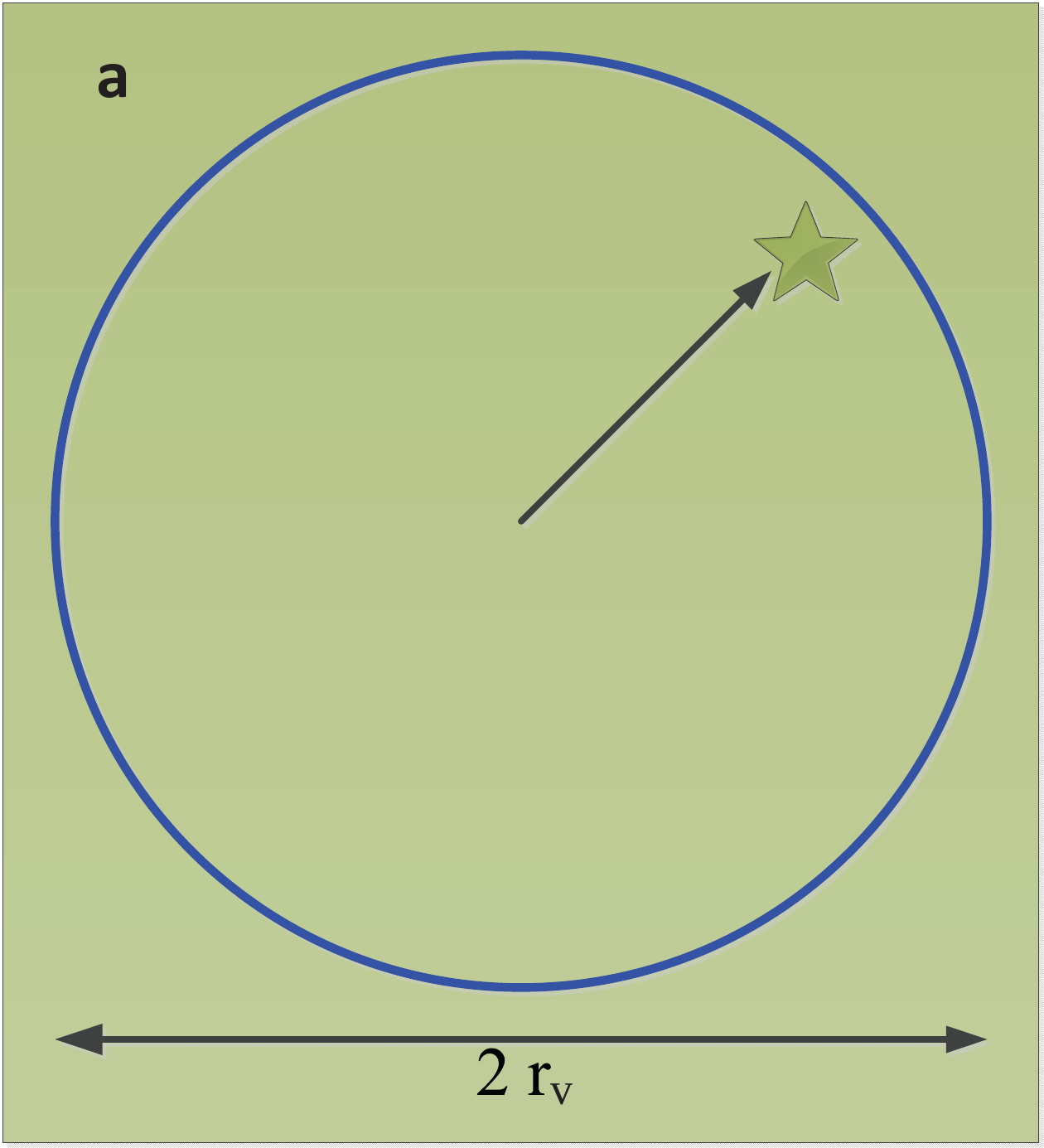}
\includegraphics[height=2in,width=2.4in]{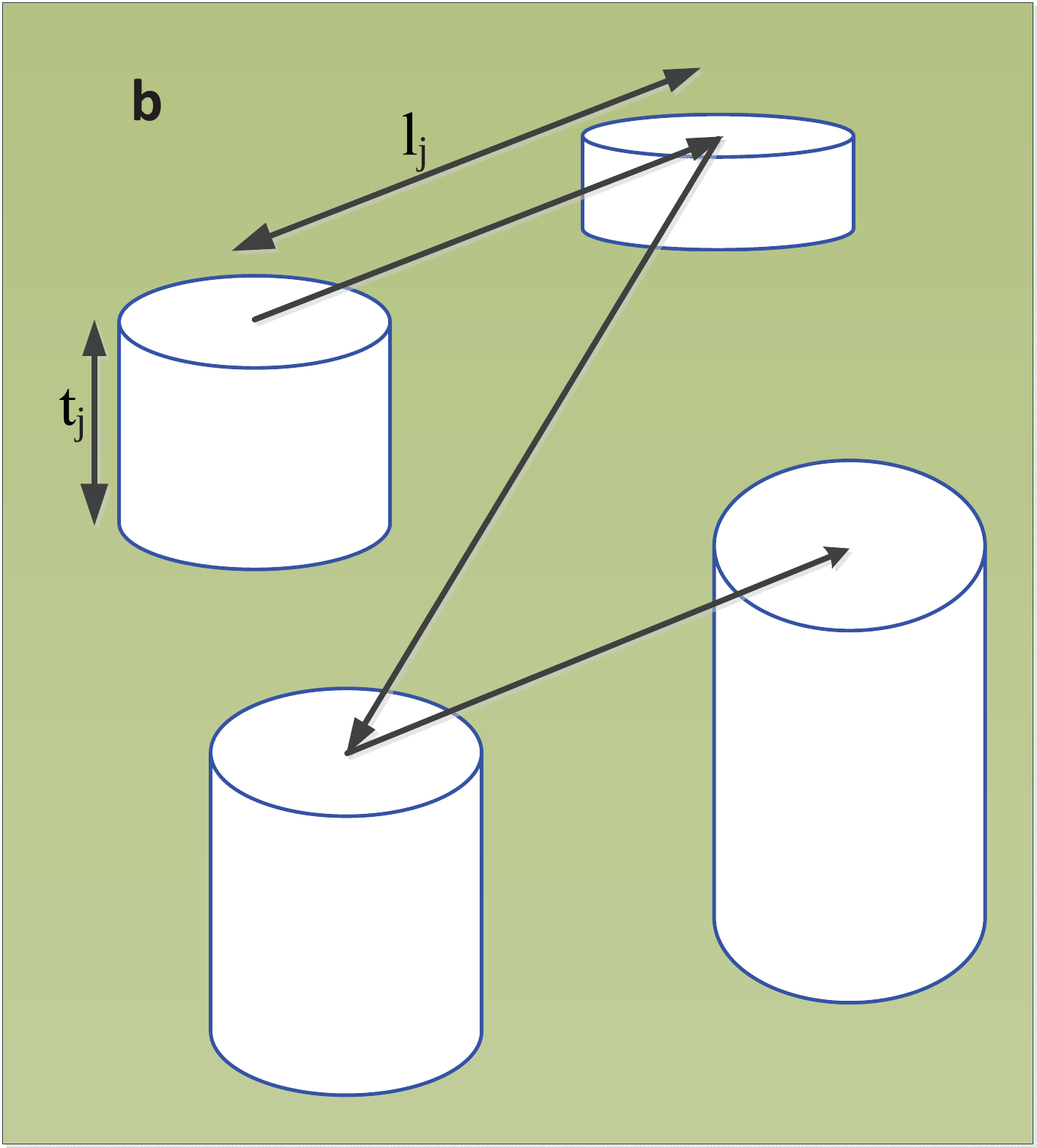}
 \bigskip

Fig. 3.1: The CTRW search strategy: (a) if the target site is
located within a visible range $r_v$, then the search agent moves on
a straight line to it; (b) if there is no target site within a
distance $r_v$, the search agent stop and wait for $t_j$ (the height
of the cylinder) and then chooses a random direction and a random
distance $l_j$ until it finds the target site.
 \end{center}

Following the similar idea of L\'{e}vy flight
\cite{Viswanathan1999,Viswanathan2000,Viswanathan2001}, we can
define the search efficiency function $\eta(\alpha, \beta)$ to be
the radio of the number of target sites visited to the total
distance traversed by the search agent as following
\begin{equation}\label{eq2.3}
\eta(\alpha, \beta)=\frac{1}{N\langle LT \rangle},
\end{equation}
where $\langle LT \rangle$ is given by
\begin{eqnarray}\label{eq2.4}
  \nonumber\langle LT \rangle &\sim& \frac{\int_{r_v}^{\lambda}l^{-\alpha}dl\int_{0}^{T}t^{-\beta}dt+
  \lambda\int_{\lambda}^{\infty}l^{-\alpha-1}dl~T \int_{T}^{\infty}t^{-\beta-1}dt}{\int_{r_v}^{\infty}l^{-\alpha-1}dl~ \int_{T}^{\infty}t^{-\beta-1}dt} \\
   &=& \frac{\alpha\beta T}{(1-\alpha)(1-\beta)}\left(\lambda^{1-\alpha}r_v^{\alpha}-r_v\right)+T \lambda^{1-\alpha}r_v^\alpha,
\end{eqnarray}
where $\lambda$ is the mean free path and $T$ is the mean wait time.
Moreover, the mean number of search $N$ between two successive
target sites for the destructive foraging case
\begin{equation}\label{eq2.5}
    N_d=\left(\frac{\lambda}{r_v}\right)^{\alpha}
\end{equation}
and for the nondestructive foraging case
\begin{equation}\label{eq2.6}
    N_n=\left(\frac{\lambda}{r_v}\right)^{{\alpha}/{2}}
\end{equation}
respectively. Substituting equations (\ref{eq2.4}) and (\ref{eq2.5})
into (\ref{eq2.3}), we find the mean efficiency $\eta(\alpha,
\beta)$ has no maximum, with lower values of $\alpha$ and $\beta$
leading to more efficient search. For the nondestructive case,
substituting (\ref{eq2.4}) and (\ref{eq2.6}) into (\ref{eq2.3}), we
find that the efficiency $\eta(\alpha, \beta)$ is optimum at
\begin{equation}\label{eq2.7}
    \alpha=1-\delta_1,~~\beta=0.5-\delta_2,
\end{equation}
where $\delta=\max(\delta_1, \delta_2)\ll 1$.

Next we test the above theoretical results with numerical
simulations. Let $r_v=1$, $T=5$, $\lambda=10, ~10^3, ~10^5, ~10^7$.
For the destructive case, the relation of the product of mean free
path and search efficiency and the parameters $\alpha$ and $\beta$
is performed in three dimensional space, shown in Figure 3.2. To
better verify the result, we just choose the 2D projection
relationship of parameters $\alpha$ and $\beta$, shown in Figure
3.3, implying that lower values of $\alpha$ and $\beta$ leading to
more efficient search, which agrees with the analytical results.

 \begin{center}
\includegraphics[height=3in,width=4.5in]{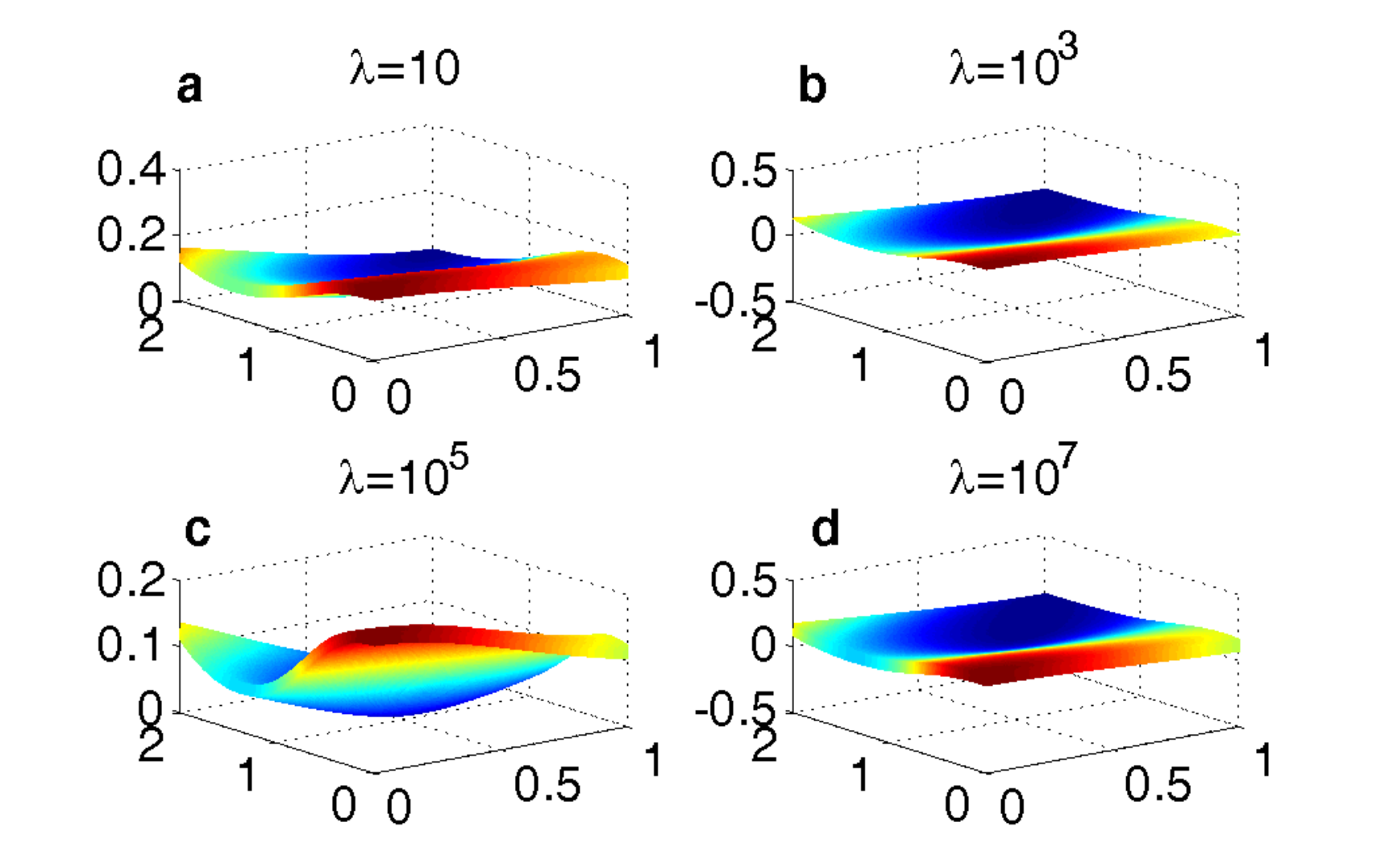}
 \bigskip

Fig. 3.2: The product of the search efficiency $\eta$ and the mean
free path $\lambda$ VS the parameter $\alpha$ and $\beta$ for
different $\lambda$: the destructive case.
 \end{center}

 \begin{center}
\includegraphics[height=3in,width=4.5in]{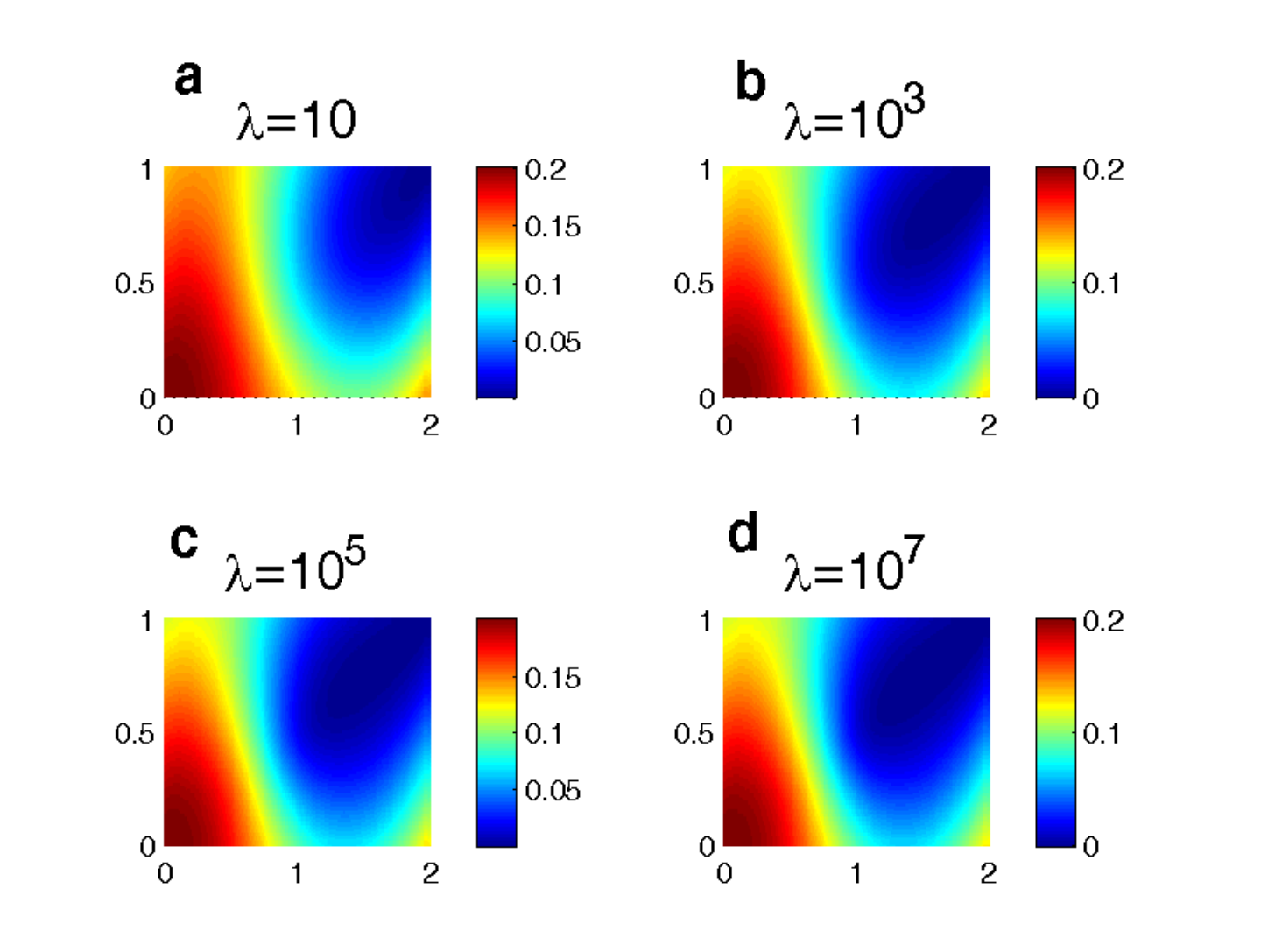}
 \bigskip

Fig. 3.3: 2D projection of Fig. 3.2 with parameter $\alpha$ and
$\beta$ for different $\lambda$.
 \end{center}

For the nondestructive case, we perform the results by using the
same parameters, shown in Figure 3.4 and Figure 3.5, from which we
find that the optimal values $\alpha=1$ and $\beta=0.5$ when
$\lambda\rightarrow\infty$. This fact shows that the analytical
results are reliable.

 \begin{center}
\includegraphics[height=3in,width=4.5in]{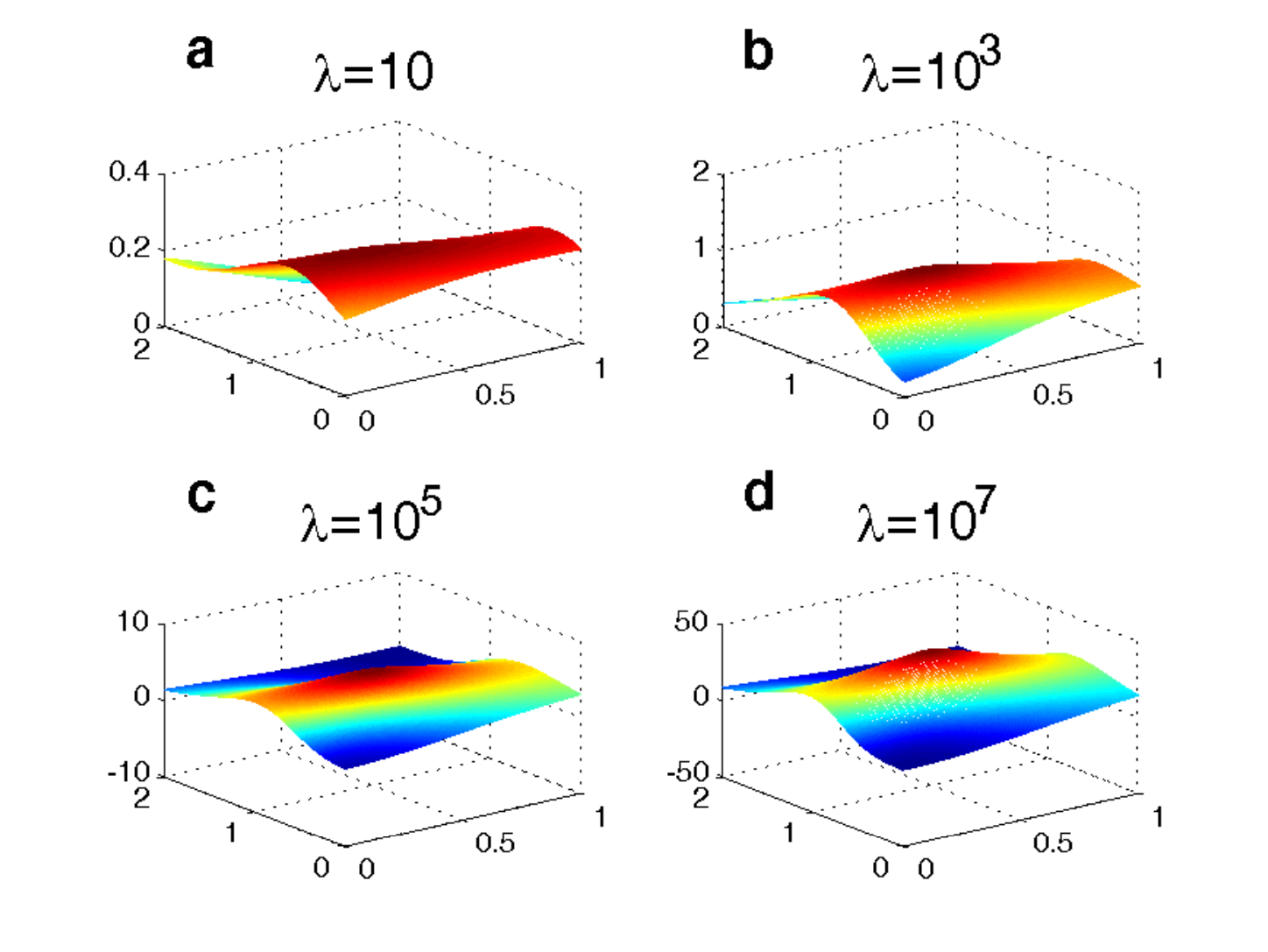}
 \bigskip

Fig. 3.4: The product of the search efficiency $\eta$ and the mean
free path $\lambda$ VS the parameter $\alpha$ and $\beta$ for
different $\lambda$: the nondestructive case.
 \end{center}

 \begin{center}
\includegraphics[height=3in,width=4.5in]{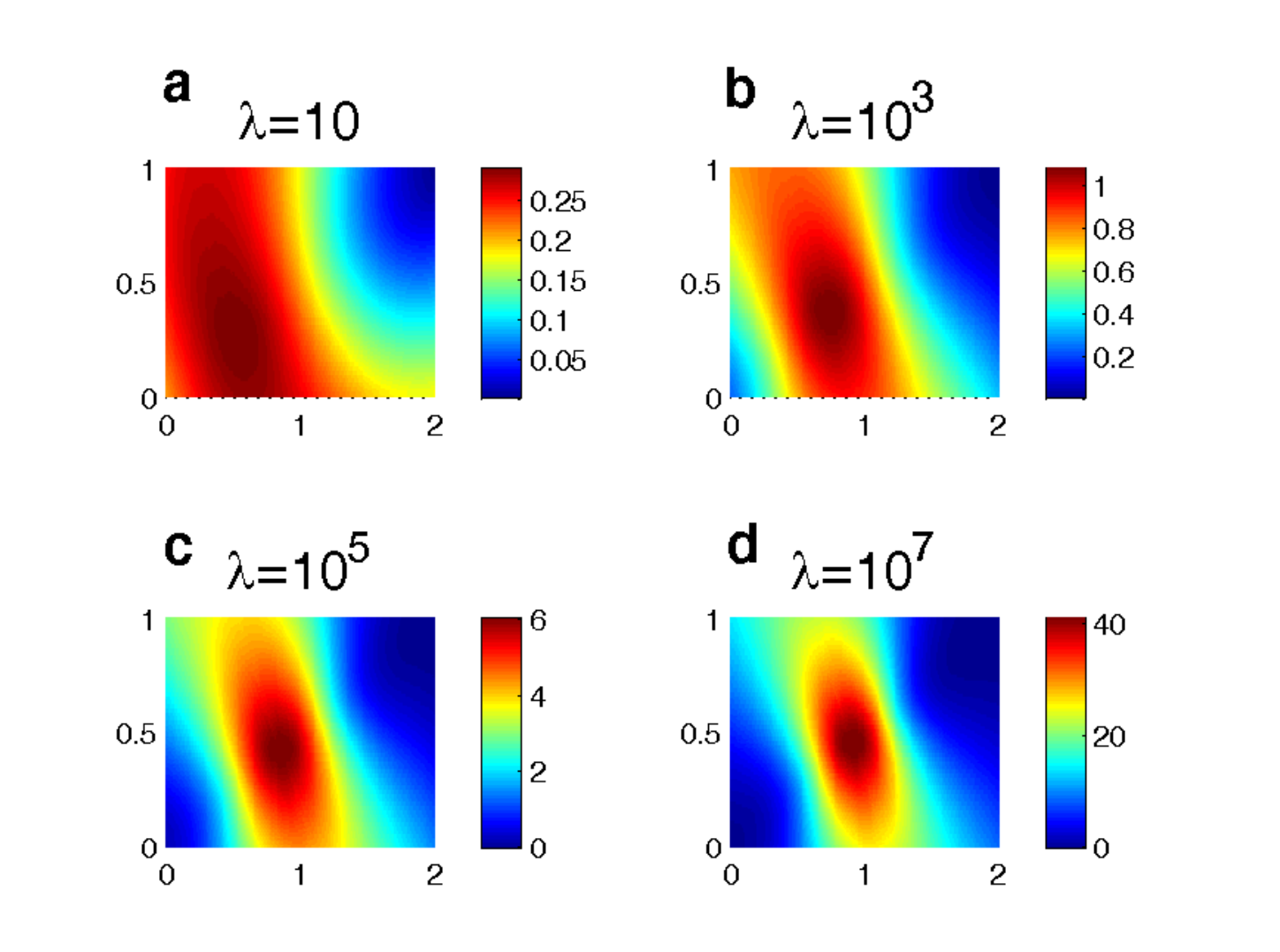}
 \bigskip

Fig. 3.5: 2D projection of Fig. 3.2 with parameter $\alpha$ and
$\beta$ for different $\lambda$.
 \end{center}

Two-dimensional path of CTRW search for $\alpha=1$ and $\beta=0.5$
is shown in Figure 3.6.

 \begin{center}
\includegraphics[height=3in,width=4.5in]{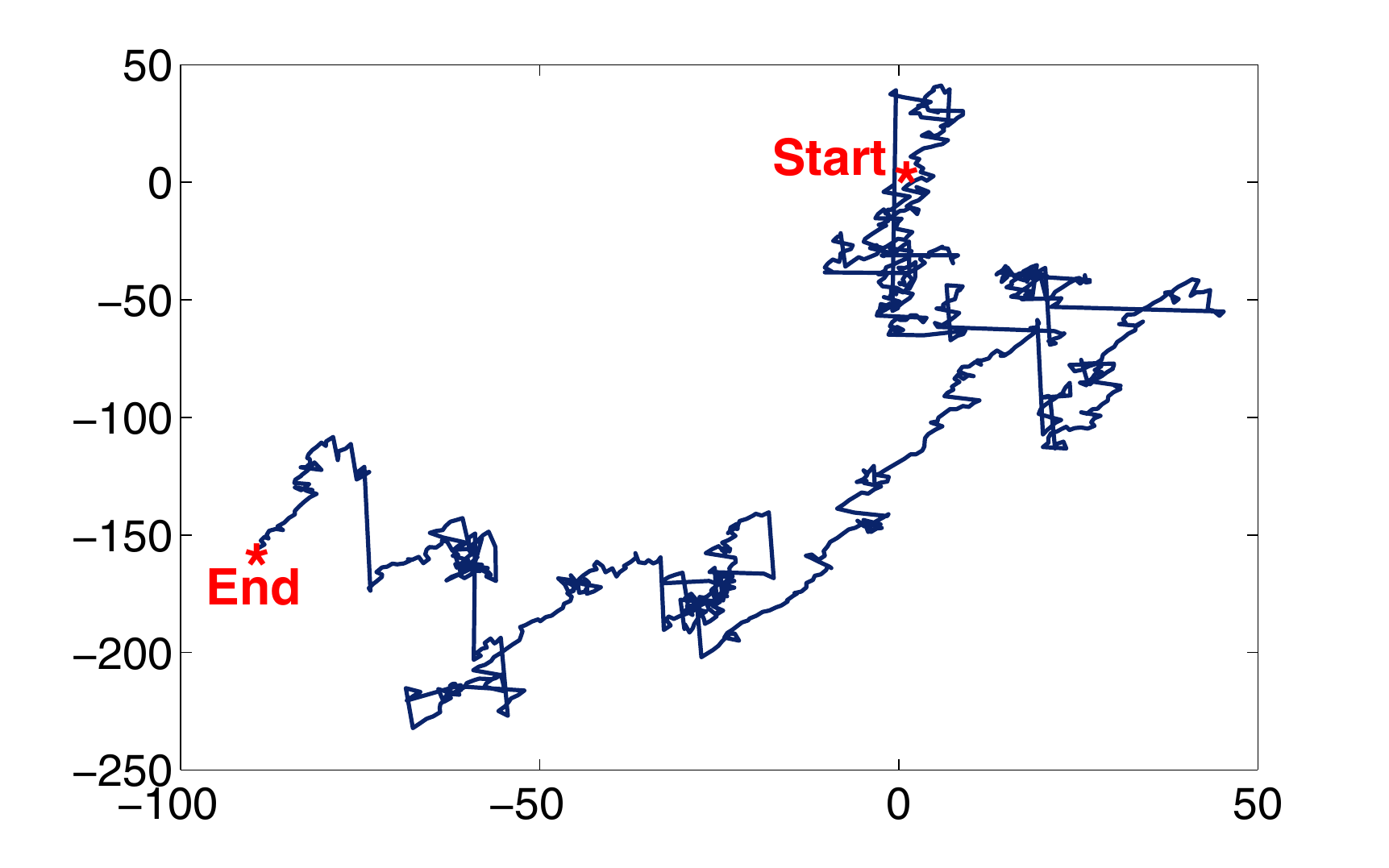}
 \bigskip

Fig. 3.6: One search path of optimal CTRW strategy in
two-dimensional space with $\alpha=1$ and $\beta=0.5$.
 \end{center}

Noting that when $\beta=0$, the search agent do not need to learn
and wait for the next step. In this case, we can ignore the waiting
time part and set $\mu=\alpha+1$, then we can cover the previous
results \cite{Viswanathan1999,Viswanathan2000,Viswanathan2001}. This
also implies that our results include the previous results
\cite{Viswanathan1999,Viswanathan2000,Viswanathan2001} as the
special case.

\section{Master equation}
\label{sec:4}

In this section, we discuss the master equation of the CTRW search
strategy in previous section based on the fractional calculus. Let
us introduce the Fourier transform for the coordinate variable $l$
and the Laplace transform for the time variable $t$:
\begin{equation}\label{eq3.1}
    W(k)=\int_{-\infty}^{\infty}e^{ikl}w(l)dl,
\end{equation}
and
\begin{equation}\label{eq3.2}
    \Psi(s)=\int_0^\infty e^{-st}\psi(t)dt
\end{equation}
respectively. Then we can get the Montroll-Weiss equation
\cite{Montroll1965,Metzler2000} in Fourier-Laplace space:
\begin{equation}\label{eq3.3}
    P(k, s)=\frac{1-\Psi(s)}{s}\frac{1}{1-W(k)\Psi(s)}.
\end{equation}
Since we assume that $w(l)$ and $\psi(t)$ are characterized by
equations (\ref{eq2.1}) and (\ref{eq2.2}), we have
\cite{Saichev1997}
\begin{equation}\label{eq3.4}
    1-W(k)\sim |k|^\alpha
\end{equation}
and
\begin{equation}\label{eq3.5}
    \frac{1}{\Psi(s)}\sim 1+s^\beta.
\end{equation}
After substituting (\ref{eq3.4}) and (\ref{eq3.5}) into
(\ref{eq3.3}), we get
\begin{equation}\label{eq3.6}
    P(k, s)=\frac{s^{\beta-1}}{s^\beta+|k|^\alpha}.
\end{equation}
On the other hand, the so-called time-space fractional order
diffusion equation is given by
\begin{equation}\label{eq3.7}
    {}_tD^\beta_\ast p(l, t)={}_lD^\alpha p(l, t),
\end{equation}
where ${}_tD^\beta_\ast$ is the time-fractional Caputo derivative of
order $\beta$, and ${}_lD^\alpha$ is the space-fractional
Riesz-Feller derivative of order $\alpha$. In fact, the Caputo
fractional derivative of $f(t)$ is defined as
\cite{Caputo1967,Podlubny1999}
\begin{equation}\label{eq3.8}
    {}_tD^\beta_\ast f(t)=\frac{1}{\Gamma(1-\alpha)} \int_{0}^{t}{(t-\tau)}^{-\alpha} \frac{d f(\tau)}{d \tau}
 d\tau,
\end{equation}
with its Laplace transform
\begin{equation}\label{eq3.9}
    \mathcal {L}\left\{{}_tD^\beta_\ast f(t); s\right\}=s^\beta F(s)-s^{\beta-1}f(0^+).
\end{equation}
The space-fractional Riesz-Feller derivative of $g(l)$ is defined as
\cite{Mainardi2001,Mainardi2010}
\begin{equation}
{}_lD^\alpha g(l)=\frac{\Gamma(1+\alpha)}{\pi}\sin\left(\frac{\alpha
\pi}{2}\right)\int_0^\infty
\frac{g(l+\xi)-2g(l)+g(l-\xi)}{\xi^{1+\alpha}}d\xi,
\end{equation}
with its Fourier transform
\begin{equation}\label{eq3.10}
    \mathcal {F}\left\{{}_lD^\alpha g(l); k\right\}=-|k|^\alpha G(k).
\end{equation}
From \cite[equation (2.17)]{Gorenflo2003}, the Laplace-Fourier
transform of fractional order equation (\ref{eq3.7}) is
\begin{equation}\label{eq3.11}
    P(k, s)=\frac{s^{\beta-1}}{s^\beta+|k|^\alpha}.
\end{equation}
which is the same as (\ref{eq3.6}). This means the proposed CTRW
search strategy obeys a time-space fractional diffusion equation.
Therefore, many interesting mathematical and physical contributions
\cite{Zeng2011,Zeng2012b,Sheng2012} on fractional dynamics can guide
us to understand the mechanism of CTRW search strategy.


\section{Concluding remarks}
\label{sec:5}

In this paper we have proposed the Continuous Time Random Walk
(CTRW) optimal search framework by assuming that both of search
length's and waiting time's distribution satisfy a power-law
function. By introducing the efficiency function, we found that the
optimum for parameters $\alpha=1$ and $\beta=0.5$ when mean free
path $\lambda$ tends to infinity. Based on fractional calculus
technique, we further derive its master equation as a time-space
fractional order diffusion equation. Thus many interesting
contributions related to fractional dynamics can guide us to
understand the mechanism of CTRW search strategy. Numerous
simulations are provided to illustrate the non-destructive and
destructive cases, which verify our analytical results are reliable.

Optimal random search is a relatively new field, and considerable
efforts are still being made and many challenges are still to be
overcome. What about the optimal search for the big data in
small-world network? How to specify the probability distributions of
search's length and waiting time in the practical applications? What
is the globally optimum search strategy for the moving target sites?
We will focus on this topic and hope this paper can stimulate wide
discussion and lead to new investigations of these challenging
problems.

\smallskip
\section*{Acknowledgements}

This work was partly supported by the National Natural Science
Foundation of China (No. 11301090, No. 11271139, No. 61104138),
Guangdong Natural Science Foundation (No. S2011040001704).




 \bigskip \smallskip

 \it

 \noindent
$^1$ School of Sciences and School of Automation Science and Engineering\\
South China University of Technology \\
Guangzhou 510640, China \\[4pt]
e-mail: zeng.cb@mail.scut.edu.cn (C. Zeng)
\hfill Received: October 18, 2013 \\[12pt]
$^2$ Mechatronics, Embedded Systems and Automation (MESA) Lab\\
School of Engineering\\
University of California, Merced\\
5200 North Lake Road, Merced, CA 95343, USA\\[4pt]
e-mail: yangquan.chen@ucmerced.edu (Y.Q. Chen)

\end{document}

%% file: fcaa_style.tex
\usepackage{graphicx}
\usepackage{epsfig}
\usepackage{amsthm}
\usepackage{amsmath}
\usepackage{latexsym}
\usepackage{amsfonts}
\usepackage{amssymb}

\newtheoremstyle{theorem}
  {15pt}          
  {15pt}  
  {\sl}  
  {\parindent}
  {\sc}  
  {. }    
  { }    
  {}     
\theoremstyle{theorem}

\newtheoremstyle{defi}
  {15pt}          
  {15pt}  
  {\rm}  
  {\parindent}     
  {\sc}  
  {. }    
  { }    
  {}     
\theoremstyle{defi}

 \def\proofname{\indent\rm P\ r\ o\ o\ f.}
 
 \def\qed{\hfill$\Box$} 
